\documentclass{article}
\usepackage{amssymb}
\usepackage{amsmath}
\usepackage{amsthm}
\usepackage{amsfonts,url}
\newtheorem{thm}{Theorem}[section]
\newcommand{\T}{\mathcal{T}}
\newcommand{\A}{\mathcal{A}}
\newcommand{\B}{\mathcal{B}}
\newcommand{\V}{\mathcal{V}}
\newcommand{\incompGraphLambda}{$G(A^{\lambda}_l,B^{\lambda}_l)$}
\newcommand{\myargmin}[2]{ 
\begin{array}{c l}
	\textrm{argmin} & #2 \\
	#1
\end{array}
}

\providecommand{\norm}[1]{\lVert#1\rVert}

\date{July 1, 2014}

\begin{document}

\title{Dynamic Geodesics in Treespace\\ via Parametric Maximum Flow}
\author{Sean Skwerer and Scott Provan}
\maketitle

\begin{abstract}
Shortest paths in treespace, which represent minimal deformations between trees, are unique and can be computed in polynomial time.
The ability to quickly compute shortest paths has enabled new approaches for statistical analysis of populations of trees and phylogenetic inference.
This paper gives a new algorithm for updating geodesic paths
when the end points are dynamic. 
Such algorithms will be especially useful when optimizing for objectives that are functions of distances from a search point to other points e.g. for finding a tree which has the minimum average distance to a collection of trees. 
Our method for updating treespace shortest paths is based on parametric sensitivity analysis of the maximum flow subproblems that are optimized when solving for a treespace geodesic.
\end{abstract}

\section{Introduction}

Evolutionary histories or hierarchical relationships are often represented graphically as phylogenetic trees.  In biology, the evolutionary history of species or operational taxonomic units (OTU's) is represented by a tree. The root of the tree corresponds to a common ancestor. Branches indicate speciation of a nearest common ancestor into two or more distinct taxa. The leafs of the tree correspond to the present species whose history is depicted by the tree.
The space of phylogenetic trees introduced by Billera, Holmes and Vogtman \cite{BHV} is a metric space in which each point corresponds to a hypothetical evolutionary history. 
Treespace has also been used in statistical analyses where populations of lungs \cite{feragen2012hierarchical} and arteries \cite{Skwerer2014} are modeled as trees.
Treepsaces are a special case of globally non-positively curved cubical complexes (also known as CAT(0) cubical complexes) \cite{BridsonHaefliger}. 
Geodesics in these spaces have applications in robotic motion planning \cite{ardila2014moving,ardila2012geodesics}.

In addition to efficient methods for computing shortest length paths \cite{Owen2011,OwenProvan}, methodologies for
other optimization problems have been developed. 
Research into methodology for treespace optimization problems includes:
\begin{enumerate}
	\item mathematical analysis of the geometry and combinatorics, as well as algorithms, for computing an average tree \cite{Miller2015}
	\item an algorithm for finding the nearest points in two convex subsets of treespace by alternating projections \cite{Bacak2012599}
	\item generalizations of the proximal point algorithm to lower-semi continuous convex functions globally non-positively curved metric spaces \cite{bacak2013proximal} and specific versions of proximal point methods which have been analyzed for functions on treespace \cite{bacak2014computing,bacak2014convex} 
	\item a generalization of principal component analysis called principal geodesics \cite{nye2011principal,nye2014algorithm}
	\item Bayesian inference for phylogenetics \cite{nye2015convergence}
\end{enumerate}

In this article we give our algorithm for dynamically updating the geodesic from a tree $T$ to a dynamic point $X$.
This method will be especially useful for accelerating optimization
routines for problems with functions of distances involving variable trees.

\section{Treespace Geodesics}\label{sec:Geodesics}
We now give an explicit description of geodesics in treespace.
Let $\T_r$ be the geometric treespace from \cite{BHV} in which each point represents a phylogenetic tree having leafs in bijection with a fixed label set $\{0,1,2,\ldots,r\}$.
Let $X \in \T_r$ be a variable point and let $T \in \T_r$ be a fixed point.
Let $\Gamma_{XT}= \{\gamma(\lambda)|0\leq \lambda \leq 1\}$ be the geodesic path from $X$ to $T$.
Let $C$ be the set of edges which are compatible in both trees, that is the union of the largest subset of $E_X$ which is compatible with every edge in $T$ and the largest subset of $E_T$ which is compatible with every edge in $X$. 

The following notation for the Euclidean norm of the lengths of a set of edges $A$ in a tree $T$ will be used frequently,
\begin{equation}
||A||_T = \sqrt{ \sum_{e \in A}{|e|_T^2}}
\end{equation}
or without the subscript when it is clear to which tree the lengths are from.

A support sequence is a pair of disjoint partitions, $A_1\cup \ldots \cup A_k =E_X\setminus C$ and $B_1\cup\ldots\cup B_k=E_T \setminus C$.
\begin{thm}\cite{OwenProvan}
	A support sequence $(\A,\B)=(A_1,B_1),\ldots,(A_k,B_k)$ corresponds to a geodesic if and only if it satisfies the following three properties:
	\begin{itemize}
		\item[\rm (P1)] For each $i>j$, $A_i$ and $B_j$ are compatible
		\item[\rm (P2)] $\frac{\norm{A_1}}{\norm{B_1}} \leq \frac{\norm{A_2}}{\norm{B_2}} \leq \ldots \leq \frac{\norm{A_{k}}}{\norm{B_{k}}}$
		\item[\rm (P3)] For each support pair $(A_i, B_i)$, there is no nontrivial partition $C_1 \cup C_2$ of $A_i$, and partition $D_1 \cup D_2$ of $B_i$, such that $C_2$ is compatible with $D_1$ and $ \frac{\norm{C_1}}{\norm{D_1}} < \frac{\norm{C_2}}{\norm{D_2}}$
	\end{itemize}
	The geodesic between $X$ and $T$ can be represented in ${\mathcal T}_r$ with legs
	\begin{displaymath}
	\Gamma_l=\left\{\begin{array}{ll}
	\left[\gamma(\lambda):\; \frac{\lambda}{1-\lambda}\leq\frac{\norm{A_1}}{\norm{B_1}}\right],
	&l=0\\[.7em]
	\left[\gamma(\lambda):\; \frac{\norm{A_i}}{\norm{B_i}}\leq\frac{\lambda}{1-\lambda}\leq\frac{\norm{A_{i+1}}}{\norm{B_{i+1}}}\right],
	&l=1,\ldots,k-1,\\[.7em]
	\left[\gamma(\lambda):\; \frac{\lambda}{1-\lambda}\geq\frac{\norm{A_k}}{\norm{B_k}}\right],
	&l=k\end{array}\right.
	\end{displaymath}
	The points on each leg $\Gamma_l$ are associated with tree $T_l$ having edge set
	
	\begin{displaymath}
	\begin{array}{rcl}
	E_l&=&B_1\cup\ldots\cup B_l\cup A_{l+1}\cup\ldots\cup A_k\cup C
	\end{array}
	\end{displaymath}
	
	Lengths of edges in $\gamma(\lambda)$ are
	
	\begin{displaymath}
	|e|_{\gamma(\lambda)}=\displaystyle\left\{\begin{array}{ll}
	\frac{(1-\lambda)\norm{A_j}-\lambda \norm{B_j}}{\norm{A_j}}|e|_X&e\in A_j\\[1em]
	\frac{\lambda \norm{B_j}-(1-\lambda)\norm{A_j}}{\norm{B_j}}|e|_{T}&e\in B_j\\[1.5em]
	(1-\lambda)|e|_X+\lambda |e|_{T}&e\in C\\
	\end{array}.\right.
	\end{displaymath}
	
	The length of $\Gamma$ is
	\begin{equation}\label{pathlength}
	d(X,T)=\bigg\Arrowvert(\norm{A_1}+\norm{B_1},\ldots,\norm{A_k}+\norm{B_k},|e_C|_{_X}-|e_C|_{_{T}})\bigg\Arrowvert
	\end{equation}
	and we call this the geodesic distance from $X$ to $T$.
\end{thm}

\section{Updating Geodesics}\label{sec:updating}
In this section we present our main results, which are sensitivity analysis of the geodesic optimality conditions and a network flow algorithm for updating the geodesic to a variable point.

\subsection{Setup and notation}
All discussion in this section takes place in squared treespace in which the coordinates of edge lengths are squared to simplify notation, and
solutions can be mapped back to the original coordinate system.
Let $X^0$ and $X^1$ be fixed trees in the same orthant,
so that the geodesic between them is a line segment. 
Let $X^\lambda = (1-\lambda)X^0+\lambda X^1$ be a variable tree on this segment.
The length of each edge in  $X^\lambda$ is $|e|_{X^\lambda}= (1-\lambda)|e|_{X^0}+\lambda |e|_{X^1}$, and the change in the length of $e$ with respect to $\lambda$ is $d_e = |e|_{X^1}-|e|_{X^0}$. 
Thus  $|e|_{X^\lambda}= |e|_{X^0}+\lambda d_e$.
Let $\Gamma^{\lambda}$ be the geodesic from $X^\lambda$ to $T$ with supports $(\A^{\lambda},\B^{\lambda}) = (A^{\lambda}_1,B^{\lambda}_1),\ldots,(A^{\lambda}_{k^{\lambda}},B^{\lambda}_{k^{\lambda}}) $.
These supports will be constant in the vistal cell $\V^\lambda$ containing $X^\lambda$. We describe conditions under which $X^\lambda$ leaves $\V^\lambda$, adn the associated updates to the supports. 

\subsection{Intersections with (P2) constraints}
The (P2) bounding inequalities for $\V^\lambda$ can be written in the form 
\begin{equation}\label{(P2)basic}
\begin{array}{c c}
\displaystyle \norm{B_{l+1}^{\lambda}}^2 \sum_{e \in A_l^{\lambda}}{|e|_{X^0}+\lambda d_e} \leq  \norm{B_l^{\lambda}}^2 \sum_{e \in A_{l+1}^{\lambda}}{|e|_{X^0}+\lambda d_e} & i = 1,\ldots, k^l-1, l = 1,\ldots, n
\end{array}
\end{equation}
Simplification yields
\begin{equation}
a_l\lambda + b_l \geq 0
\end{equation}
where 
\begin{displaymath}
a_l = \norm{B_{l}^{\lambda}}^2 \sum_{e \in A_{l+1}^{\lambda}}{d_e}-\norm{B_{l+1}^{\lambda}}^2 \sum_{e \in A_{l}^{\lambda}}{d_e},
\end{displaymath} and
\begin{displaymath}
b_l = \norm{B_{l}^{\lambda}}^2 \sum_{e \in A_{l+1}^{\lambda}}{|e|_{X^0}}-\norm{B_{l+1}^{\lambda}}^2 \sum_{e \in A_{l}^{\lambda}}{|e|_{X^0}}.
\end{displaymath}
There are several cases for solutions, $\lambda_l = -b_l/a_l$, and each signifies a distinct situation.
The case $\lambda=0$ implies $X^0$ is on a (P2) boundary of $\V^\lambda$.
Any positive solution, $0<\lambda \leq 1$ corresponds to a point along the geodesic segment at which the segment intersects a (P2) constraint.  
Finding $\lambda > 1$ signifies an intersection with a (P2) boundary beyond $X^1$, and
if a solution is negative then there is an intersection with the geodesic in the opposite direction to $d$.
The first (P2) inequalities to be violated in moving along the geodesic segment are 
\begin{equation}
\displaystyle \myargmin{l: \lambda_{l}>0}{ \left\{ \lambda_l\right\}}
\end{equation}
If a (P2) constraint is satisfied at equality, the corresponding supports may be combined to make a new support satisfying (P2) at strict inequality, and still satisfying (P1) and (P3). Combining the current flow values into this new support pair will provide a warm start for subsequently tracking intersections with (P3) constraints, and further (P2) violations.

\subsection{Intersections with (P3) constraints}

From \cite[Prop. 3.3]{Miller2015} inequality constraints for (P3) are
\begin{equation}\label{(P3)basic}
\begin{array}{l}
\norm{B_l^{\lambda} \setminus J}^2 \sum_{e \in A_l^{\lambda} \setminus I} |e|_{X^0}+\lambda d_e - \norm{J} ^2\sum_{e \in I} |e|_{X^0}+\lambda d_e \geq 0 \\ \;for\;all\; i = 1,\ldots,k \; and\; subsets \;
I\subset A_l^{\lambda},\;J\subset B_l^{\lambda} \; such \; that \; I \cup J \; is \; \\compatible.
\end{array}
\end{equation}
Determining whether or not a support pair satisfies (P3) can be restated as the following extension problem:

\noindent{\bf Extension Problem}

\noindent{\bf Given:} Sets $A \subseteq E_X$, and $B \subseteq E_{T}$

\noindent{\bf Question:} Does there exist a partition $C_1  \cup C_2$ of $A$ and a partition of $D_1 \cup D_2$ of $B$ such that
\begin{itemize}
	\item[(i)] $C_2 \cup D_1$ corresponds to an independent set in $G(A,B)$,
	\item[(ii)] $\frac{ \norm{C_1}}{\norm{D_1}} \leq \frac{\norm{C_2}}{\norm{D_2}}$
\end{itemize}
The extension problem can be reformulated from a maximum independent set problem to a maximum flow problem. 
Each support pair $(A_l^{\lambda},B_l^{\lambda})$ has a corresponding incompatibility graph as defined in \cite[Sec. 3]{OwenProvan}.
The vertex weights of the incompatibility graph at a point along the geodesic segment can be parameterized in terms of $\lambda$ as
\begin{equation}
w_e^\lambda = 
\left \{ \begin{array}{c c}
\frac{|e|_{X^0}+\lambda d_e}{\sum_{e' \in A_l^{\lambda}} |e'|_{X^0}+\lambda d_e'} & e \in A_l \\
\\
\frac{|e|_{T}}{\sum_{e'\in B^{\lambda}_l}{|e'|_{T}}} & e \in B_l 
\end{array} \right .
\end{equation}
Although $w_e^\lambda$ is a non-linear function of $\lambda$,
matters are simplified by rescaling the lengths of edges to have sum 1 within each support pair separately. 
The approach is to complete the parametric analysis of
that extension problem and then map back 
to find $\lambda$ for the original weights.

Let $V^0$ and $V^1$ be formed
by rescaling the lengths of edges in $X^0$ and $X^1$ to sum to 1. 
Suppose that some (P3) constraint defined by $\sum b_e |e|_X \geq 0 $ is satisfied at equality by $\tilde{\lambda}$, that is $\sum b_e((1-\tilde{\lambda})|e|_{V^0}+\tilde{\lambda} |e|_{V^1})=0$.  The following gives a transformation between the solution when the weights in each support pair are scaled to have sum 1, and before scaling. 

\begin{align}
\sum b_e((1-\tilde{\lambda})|e|_{V^0}+\tilde{\lambda} |e|_{V^1})=0\\
\sum b_e\left((1-\tilde{\lambda})\frac{|e|_{X^0}}{\sum_{e' \in A_l}{|e'|_{X^0}}}+\tilde{\lambda} \frac{|e|_{X^1}}{\sum_{e' \in A_l}{|e'|_{X^1}}}\right)=0\\
\sum b_e\left(c_0 |e|_{X^0}+c_1 |e|_{X^1}\right)=0\\
\sum b_e\left((1-\lambda)|e|_{X^0}+\lambda |e|_{X^1}\right)=0
\end{align}
Thus $\sum b_e((1-\lambda)|e|_{X^0}+\lambda |e|_{X^1})=0$ is
satisfied by 
\begin{align}
\lambda = \frac{c_1}{c_0+c_1} = \frac{\tilde{\lambda}/\sum_{e \in A_l}{|e|_{X^1}}}{(1-\tilde{\lambda})/\sum_{e' \in A_l}{|e'|_{X^0}}+\tilde{\lambda}/\sum_{e' \in A_l}{|e'|_{X^1}}}.
\end{align}

Assuming the weights of edges in each support pair are already scaled to have 
sum 1, the weights are parameterized as a linear function as
\begin{equation}
\tilde{w}_e^{\tilde{\lambda}} = 
\left \{ \begin{array}{c c}
|e|_{A_l}+\tilde{\lambda} \tilde{d}_e  & e \in A_l \\
\\
|e|_{B_l} & e \in B_l 
\end{array} \right .
\end{equation}
where $\tilde{d}_e = \frac{|e|_{X^1}}{\sum_{e' \in A_l}{|e'|_{X^1}}}-\frac{|e|_{X^\lambda}}{\sum_{e' \in A_l}{|e'|_{X^0}}}$ is the change in capacity for the arc from the source to node $e \in A_l^{\lambda}$.

Parametric analysis of the extension problem will yield a method for 
updating the objective function along the segment from $X^0$ to $X^1$. 
In an incompatibility graph
the capacity for an arc from the source $s$ to a node $e$ in $A_l$ is
$c_{se} = w_e$, for arcs from a node $f$ in $B_l$ to the terminal node $t$
the capacity is $c_{ft}=w_f$ (fixed), and the capacity for an arc
from a node in $A_l$ to an incompatible node in $B_l$
is infinity. 
Assume an initial flow for the directed graph \incompGraphLambda\, is calculated for $\lambda = 0$.
For arc $(e,e')$, variable flow is $z_{ee'}$ and the flow for $\lambda = 0$ is $z_{ee'}^0$.
Residual capacity for arc $(e,e')$ is $r_{ee'}= c_{ee'}-z_{ee'}+z_{e'e}$.
Recall, that $\tilde{d}_e = |e|_{V^1}-|e|_{V^0}$ is the change in capacity for the arc from 
the source to node $e \in A_l^{\lambda}$ from $X^0$ to $X^\lambda$.
As $\lambda$ increases the flow may become infeasible because an arc capacity has decreased, or
an augmenting path may exist because the arc capacity has increased. 
The net change in the total arc capacity is zero because $\sum_e \tilde{d}_e = \sum_e (|e|_{V^1}-|e|_{V^0}) = 0$, and thus
the total increase in arc capacity must equal the total decrease in arc capacity.
Therefore the maximum flow will be inhibited not by a change in total capacity, but rather by a bottleneck preventing
a balance of flow as $\lambda$ increases. 
To balance the flow, excess flow from arcs with decreasing capacities must shift to arcs with increasing capacities.
The key is to identify directed cycles in the residual graph oriented along arcs with increasing capacities
and against arcs with decreasing capacities.

If $d_e>0$ then $e$ is a ``supply'' node and if $d_e < 0$ then $e$ is a ``demand" node. 
A (P3) constraint is violated precisely at the smallest positive $\lambda$ such that balancing supply and demand is not possible.
An augmenting path is a path in the residual network
from a supply node to a demand node. If there is an augmenting path
from each supply node to each demand node then it is possible to maintain
a feasible flow for some $\lambda >0$ by pushing flow along augmenting paths. 
In pushing flow along a set of augmenting paths, where $P_e$ is the
augmenting path for supply node $e$, the residual capacity along 
arc $(e',e'')$ is 
\begin{align}
r_{e'e''}= r^0_{e'e''}+\lambda \left( \sum_{e:(e'',e')\in P_e} d_e-\sum_{e:(e',e'')\in P_e} d_e\right) 
\end{align}
For a set of augmenting paths an arc is a bottleneck at $\lambda$ if it has
has zero residual capacity at $\lambda$.
Once a bottleneck is reached at least one augmenting path
is no longer valid.
Thus a given set of augmenting paths
cannot feasibly balance the total flow at
the smallest positive value $\lambda^*$ which has a bottleneck arc.

Each supply node with flow blocked by a bottleneck needs
a new augmenting path. If such a path cannot be found
then supply and demand cannot be balanced for $\lambda > \lambda^*$.
This signifies that $X^{\lambda^*}$ is on a (P3) boundary of its vistal cell. Thus support pair $(A_l,B_l)$
could be partitioned into a support pairs $(C_1,D_1)$ and $(C_2,D_2)$ (or even
into a sequence of support pairs as described in \cite[Lem. 3.23]{Miller2015}) to create
a valid support for the same geodesic from $X^\lambda$ to $T$.
If a supply node does not have any augmenting path, then increasing $\lambda$
will result in excess flow capacity which cannot be utilized to push flow from $s$ to $t$.

A minimum cover for $\lambda >\lambda^*$ can be constructed; and in what follows
``the minimum cover" refers to the one which is being constructed.
If $e$ does not have an augmenting path, then increasing $\lambda$ will cause
the residual capacity $r_{se}$ to become positive in a maximum flow.
Therefore $e$ cannot be part of the minimum cover. Consequently, to cover
the edges adjacent to $e$, every node adjacent
to $e$ in $B_l$ must be in the minimum cover. 
Supply nodes which have augmenting paths will be in the minimum cover
unless all of their adjacent arcs are adjacent to nodes in $B_l$ which are already
in the minimum cover. 
In summary the minimum cover is $C_1 \cup D_2$ where
$D_2$ is comprised of elements in $B_l$ which are adjacent to elements of $A_l$ which
do not have augmenting paths and 
$C_1$ is comprised of elements
in $A_l$ which do not have augmenting paths or which
have all their adjacent arcs covered by elements from $B_l$. 
The new support sequence is formed by replacing $(A_l, B_l)$ with 
$(C_1,D_1),(C_2,D_2)$ where $C_2 = A_l \setminus C_1$ and $D_1 = B_l \setminus D_2$.

Standard net flow techniques can be used to find a feasible flow from
supply nodes to demand nodes. If no feasible flow exists, then 
the cut from the Supply-Demand Theorem can be used
to construct a minimum cover for the solution to the (P3) extension problem.

	There are many choices for how to find an augmenting path $P_e$ for supply node $e$.
	One method is to find a minimum spanning tree of the residual network for each supply node. Minimum spanning tree algorithms vary in computational cost, for example
	Prim's Algorithm has complexity $O(r^2)$. Once a minimum spanning tree is established
	for supply node $e$ removing bottleneck arcs one at a time will only 
	require adjusting at most one arc to create a new minimum spanning tree.
	If the demand node adjacent to 
	the bottleneck arc has no other incoming arcs
	with positive residual capacity then
	it can no longer be reached from any supply node.
	If there are any incoming arcs adjacent to that demand node with
	positive residual capacity, then adding one
	which is adjacent to the element of $B_l$ with the shortest distance from $e$
	will create a minimum spanning tree for $e$. Therefore 
	updating the minimum spanning tree for edge $e$ has cost at most equal
	to the number of elements in $B_l$ and the total cost of updating the minimum 
	spanning trees of all supply nodes is $O(r^2)$. 


\noindent {\bf (P3) Intersection Algorithm}\\
\noindent{\bfseries initialize} $\Lambda = 0$, $r_{ee'} = c_{ee'}$, $z_{ee'} = 0$\\
\noindent \hspace*{1 cm} find a maximum flow in $G(A^{\lambda,i}_l,B^{\lambda,i}_l)$\\
\noindent{\bfseries while} $\Lambda < 1$\\
\noindent \hspace*{1 cm} {\bfseries do} \\
\noindent \hspace*{1 cm} find a feasible flow in the residual network \\
\noindent \hspace*{1 cm} {\bfseries if} supplies and demands are infeasible \\
\noindent \hspace*{2 cm} halt, a (P3) boundary intersection \\
\noindent \hspace*{1 cm} {\bfseries endif}\\
\noindent \hspace*{1 cm} augment flow until some residual capacity reaches zero\\
\noindent \hspace*{1 cm} calculate smallest $\lambda^* > \Lambda$ with a bottleneck arc \\
\noindent \hspace*{1 cm} using the residual capacities:\\
\noindent \hspace*{1 cm} $r_{e'e''}^{\lambda^*}= r^\Lambda_{e'e''}+(\lambda^*-\Lambda) \left( \sum_{e:(e'',e')\in P_e} d_e-\sum_{e:(e',e'')\in P_e} d_e\right) $\\
\noindent \hspace*{1 cm} $\Lambda = \lambda^*$\\
\noindent{\bfseries endwhile}

	Initialization, setting residual capacities, and finding a maximum flow has
	computational complexity $O(r^3)$.
	The (P3) Intersection Algorithm halts in fewer than $r^2$ iterations of the while loop.
	Once a bottleneck arc has residual capacity zero it will stay zero.
	If all of the bottleneck arcs, of which there are fewer than $r^2$, reach residual capacity zero then there are no augmenting paths. This is more
	than sufficient to cause at least one supply node to not have any 
	augmenting path.
	Calculating $\lambda^*$ requires updating
	$ \sum_{e:(e'',e')\in P_e} d_e-\sum_{e:(e',e'')\in P_e} d_e $
	when the augmenting path for $e$ changes. Updating these terms
	when necessary requires identifying when the augmenting path for $e$
	changes, removing $d_e$ for each arc on the old augmenting path for $e$,
	and including $d_e$ for each arc on the new augmenting path for $e$.
	In the worst case the augmenting paths for all supply nodes will change
	so the total cost of updating the rate of change in flow along every arc is $O(r^2)$.
	Therefore the (P3) Updating Algorithm will either reach $\lambda=1$ 
	or halt in fewer than $r^2$ iterations of the while loop.
	The total complexity for the (P3) Intersection Algorithm using
	minimum spanning tree updating to find augmenting paths is  $O(r^4)$.

\section{Conclusions and future directions}
Optimization for non-positively cubical complexes is an area which has
received attention due to applications these spaces in modeling trees, in phylogenetics and in robotic motion planning.
The challenges of solving optimization problems quickly must be overcome.
Our methodology for updating geodesics has the potential to accelerate optimization methods for problems with distances involving dynamic trees.
This includes computing averages, principal geodesics, finding nearest points in two convex sets and Bayesian inference for phylogenetic trees. 

The incorporation of our method for updating geodesics to accelerate
existing methods is an important topic of further research.
Several important questions remain.
Some of the optimization methods may only make small adjustments, however others, such as proximal point algorithms, can make drastic changes to the search trees. It remains to be determined, on the basis of specific problems, whether it will be more efficient to update or to start from scratch.
The question is open, if and how our method can be generalized to non-positively curved cubical complexes other than treespace.

\bibliographystyle{acm}
\bibliography{thesis,citations,library,trees,treespaceApplications,FKS}

\end{document}